\documentclass[12pt]{amsart}
\usepackage{amsmath,amsthm,amssymb,color}

\newtheorem{theorem}{Theorem}[section]
\newtheorem{thm}{Theorem}[section]
\newtheorem{lemma}[theorem]{Lemma}
\newtheorem{cor}[theorem]{Corollary}
\theoremstyle{remark}
\newtheorem{remark}[theorem]{Remark}
\theoremstyle{definition}
\newtheorem{definition}[theorem]{Definition}
\theoremstyle{proposition}
\newtheorem{proposition}[theorem]{Proposition}
\numberwithin{equation}{section}
\allowdisplaybreaks

\begin{document}
\title{Abelian Functions for Cyclic Trigonal Curves of Genus Four}
\author{S. Baldwin}
\address{Department of Mathematics, Imperial College, London, UK SW7 2AZ}
\email{sadie.baldwin@imperial.ac.uk}
\author{J. C. Eilbeck}
\address{The Maxwell Institute and Department of Mathematics, 
Heriot-Watt University, Edinburgh, UK EH14 4AS}
\email{J.C.Eilbeck@hw.ac.uk}
\author{J. Gibbons}
\address{Department of Mathematics, Imperial College, London, UK SW7 2AZ}
\email{j.gibbons@imperial.ac.uk}
\author{Y. \^Onishi}
\address{Faculty of Humanities and Social Sciences,
Iwate University, Ueda 3-18-34, Morioka 020-8550, Japan}
\email{onishi@iwate-u.ac.jp}
\begin{abstract}
  We discuss the theory of generalized Weierstrass $\sigma$ and $\wp$
  functions defined on a trigonal curve of genus four, following
  earlier work on the genus three case.  The specific example of the
  ``purely trigonal'' (or ``cyclic trigonal'') curve
  $y^3=x^5+\lambda_4 x^4 +\lambda_3 x^3+\lambda_2 x^2 +\lambda_1
  x+\lambda_0$ is discussed in detail, including a list of some of the
  associated partial differential equations satisfied by the $\wp$
  functions, and the derivation of an addition formulae.
\end{abstract}

\maketitle
\section{Introduction}
Over the last few years there has been increasing interest in
explicit descriptions of Abelian functions of algebraic curves.
The beginnings of this theory go back to Weierstrass's theory of elliptic
functions, and we will take this as our model. The key results of this are
simply stated: let $\sigma(u)$ and $\wp(u)$ be the standard functions in
Weierstrass elliptic function theory.  They satisfy the well-known formulae
\begin{eqnarray}
\wp(u) &=& - \frac{d^2}{du^2}\ln \sigma (u),\\
(\wp')^2&=&4\wp^3-g_2\wp-g_3,\\
\wp''&=&6\wp^2-\tfrac12 g_2 \label{WP}.
\end{eqnarray}

The $\sigma$-function has a power series expansion, originally due to
Weierstrass:
\begin{equation}
  \sigma(z) = 
  u-{\tfrac {1}{240}}\,g_{{2}}{u}^{5}- {\tfrac {1}{840}}\,g_{{3}}{u}^{7}- {
    \tfrac {1}{161280}}\,{g_{{2}}}^{2}{u}^{9}- {\tfrac {1}{2217600}}\,g_{{2}}
  g_{{3}}{u}^{11}+\dots
\end{equation}
Here the coefficients satisfy a linear recurrence, which is
given in \cite{AS} (see also \cite{ee00}).

The $\sigma$-function for an elliptic curve also satisfies a two-term 
addition formula which is a key result of the theory:
\begin{equation}
   -\frac{\sigma(u+v)\sigma(u-v)}{\sigma(u)^2\sigma(v)^2}=\wp(u)-\wp(v).
   \label{elliptic2term}
\end{equation}
Taking logarithmic derivatives, one may obtain the standard algebraic
addition formula for $\wp(u)$.

Klein's generalisation of this theory is described in the classical
monograph of Baker \cite{ba97}; in particular, in \cite{ba07} he
worked out in detail the formulae corresponding to these for a genus 2
hyperelliptic curve. More recently such addition formulae for other
curves have been derived; for example \cite{bel97}, where it was shown
that for a hyperelliptic curve of arbitrary genus, the right hand side
of the two-term addition formula for $\sigma$ generalises to an
expression written in terms of a Pfaffian.

The corresponding theory for trigonal curves is more complex.  Few
examples have been worked out explicitly: the $\sigma$-function
realization of the Abelian functions of a trigonal curve was developed
in some detail in \cite{bel00} and \cite{eel00} and some of the
present authors \cite{bg06,eemop06,on05},  have studied
specific examples.  These examples have concentrated on
cyclic trigonal curves.  This sub-set of curves possesses symmetry
properties which simplify the calculations considerably, thus
providing a useful starting point to investigate the properties of
Abelian functions of trigonal curves in more detail.

There is a natural splitting of $(3,s)$ cyclic trigonal curves, with $
s\ne 0 \mod 3$, into those with $s=3n+1$ and those with $s=3n+2.$ The
simplest non-trivial case of this first class, $s=3n+1,$ is the 
$(3,4)$ curve of genus 3.  The two papers \cite{eemop06,on05} considered the
general case for this curve and then looked in greater detail at the
corresponding cyclic curve.  For the cyclic curve, they found not only
an explicit two-term addition formula, which has the same left-hand
side as that of the Weierstrass' addition formula \eqref{elliptic2term},
but also a three-term addition formula.  This three-term addition
formula reflects the natural symmetries of the cyclic curve.  It gives
an expansion of
\begin{equation}
  \frac{\sigma(u+v)\sigma(u+[\zeta]v)\sigma(u+[\zeta^2]v)}
  {\sigma(u)^3\sigma(v)^3}, \qquad
  \zeta = \exp\left(\frac{2 \pi i}{3}\right) \nonumber
\end{equation} 
in terms of Abelian functions of $u$ and $v$. A key tool in
calculating these addition formulae was the evaluation of two bases of
linearly independent Abelian functions.  The first of these bases has
poles of order at most two on the $\Theta$-divisor, where $\sigma$ has
a simple zero, and the second basis has poles of order at most three
on the $\Theta$-divisor.

The next case to consider is the simplest non-trivial curve of the
second class: a $(3,5)$ trigonal curve with genus 4.  We should stress
here that the $(3,4)$ and $(3,5)$ cases have some important
differences.  Most immediate is the fact that for the $(3,4)$ curve,
$\sigma$ is an odd function of $u$, while for the $(3,5)$ case it is
even. More generally for a $(3,s)$ curve where $s$ is not divisible by
$3$, $\sigma$ has odd or even parity for $s$ respectively even or odd.
This change of parity has consequences in subsequent formulae, most
clearly in the two-term addition formulae - the right hand side is
respectively antisymmetric or symmetric between $u$ and $v$ in the
$(3,4)$ and $(3,5)$ cases.

The cyclic $(3,5)$ curve was considered in \cite{bg06} and is the
topic of this paper. A series expansion of the $\sigma$-function for
the $(3,5)$ curve was found by two of the present authors in
\cite{bg06}, and that series, extended to higher order, plays a
crucial role in some of the proofs below. They also found explicit
formulae for a basis of differentials and for the Jacobi inversion
formula for the curve. A basis of Abelian functions with second order
poles has been constructed, leading to sets of differential equations
satisfied by the fundamental Abelian functions of the curve; these are
of three kinds, those expressible as a fourth-order quasilinear PDE
for $\sigma$, and those expressible as third-order PDE for $\sigma$,
either quasilinear or else of second degree.

These results have also enabled us to find an explicit two-term
addition theorem. Further, it has been possible to demonstrate that a
three-term addition theorem exists for this curve and indeed for
a cyclic trigonal curve of any genus. However, in contrast to the
$(3,4)$ case, a basis of Abelian functions with third order poles has
not been explicitly constructed, so no explicit form for the 3-term
addition formula can yet be found.

As commented in \cite{eemop06}, our study is far from complete.  One
problem still to be considered is an explicit recursive construction
of the $\sigma$-series generalizing the one given by Weierstrass. For
hyperelliptic curves of genus two this was done in \cite{bl02}.
Another goal is to understand the algebraic structure of the two kinds
of addition theorem developed here, in order to generalize results to
higher genera, as \cite{bel97} did for hyperelliptic curves.
Buchstaber and Leykin have described some progress on these problems
for general curves in \cite{bl05}.

Section 2 of this paper reviews the relevant theory of holomorphic
differentials on a trigonal, and particularly a cyclic trigonal curve;
in section 3 the $\sigma$-function is introduced and its properties
discussed - following Klein, our discussion uses the fundamental
bidifferential, which plays a major role in what follows.  In section
4 we introduce the fundamental Abelian functions $\wp_{ij}$ on the
curve as well as some important differential expressions in these,
denoted as $Q_{ijk\ell}$.  We then expand the Klein
bidifferential in two ways to obtain the Jacobi inversion formula and
some simple examples of the two kinds of PDE satisfied by the
$\wp_{ij}$. Section 5 deals with the vector space $\varGamma(J,
\mathcal{O}(2\Theta^{[2]}))$, of functions with at most double poles
on the $\Theta$-divisor, where $\sigma(u)$ vanishes. We construct an
explicit basis for this space.  In section 6 we describe the Taylor
expansion of $\sigma(u)$ near the origin, first given in \cite{bg06}.
This is then used, with the results of section 5, to obtain the PDEs
satisfied by the Abelian functions, of which the known examples are
listed in section 7.  Finally sections 8 and 9 are concerned
respectively with the two- and three-term addition formulae satisfied
by $\sigma(u)$ - the results of section 5 enable the two-term formula
in section 8 to be given explicitly.

\section{Trigonal curves of degree five} 
\setcounter{equation}{0}

We define a general trigonal curve with the unique branch point
$\infty$ at infinity by $g(x,y)=0,$ where
\begin{align}
\begin{split}
  g(x,y) =
  y^3 + &(\mu_2 x + \mu_5)y^2 + (\mu_1 x^3+\mu_4 x^2 + \mu_7 x + \mu_{10})y \\
  & - (x^5 + \mu_3x^4 + \mu_6x^3 + \mu_9x^2 + \mu_{12} x +\mu_{15})
  \quad \text{($\mu_j$ are constants).}
\end{split}\label{trig35}
\end{align}
This curve is of genus 4, if it is non-singular.  The topic of this
paper is a special subclass of this family, the curves $C$ of the form
$f(x,y)=0$, where
\begin{align}
\begin{split}
   f(x,y) =
   y^3  - (x^5 + \lambda_4 x^4 + \lambda_3 x^3 + \lambda_2 x^2 + \lambda_1
    x +\lambda_0),
    \quad \text{($\lambda_j$ are constants).}
\end{split}\label{cyclic35}
\end{align}
A curve of this form is called a {\em cyclic} trigonal 
curve, \cite{Accola} trigonal 
curve; such curves are invariant under the cyclic symmetry
$$
[\zeta]: \, (x,y)\rightarrow (x,\zeta y),
$$ 
where $\zeta$ is a cube root of unity. In \cite{on05}  and \cite{eemop06}
these curves are called {\em purely} trigonal, to distinguish them from other
trigonal curves invariant under a cyclic group.  This symmetry plays a
significant role in what follows.  All objects of the theory must
transform simply under this group action. The argument we use here is
closely modelled on that of \cite{eemop06} which studied the analogous
problem for a cyclic $(3,4)$ curve.

We consider the set of differentials $\omega_1, \ldots, \omega_4$ where 
\begin{equation}
   \omega_1=\frac{\mathrm{d}x}{\frac{\partial}{\partial y}f(x,y)},  \quad
   \omega_2=\frac{x\mathrm{d}x}{\frac{\partial}{\partial y}f(x,y)},  \quad
   \omega_3=\frac{y\mathrm{d}x}{\frac{\partial}{\partial y}f(x,y)},  \quad
   \omega_4=\frac{x^2\mathrm{d}x}{\frac{\partial}{\partial y}f(x,y)}.
   \label{eq1.2} 
\end{equation}
This is a basis of the space of differentials of the 1st kind on $C$.
We denote the vector consisting of the forms (\ref{eq1.2}) by
\begin{equation}
  \omega=(\omega_1,\omega_2,\omega_3,\omega_4)
  \label{eq1.2.5}
\end{equation}
From the general theory, we know that for four variable points
$(x_1,y_1)$, $(x_2,y_2)$, $(x_3,y_3)$, $(x_4,y_4)$ on $C$, the sum of
integrals from $\infty$ to these four points
\begin{equation}
\begin{aligned}
   u &=(u_1,u_2,u_3,u_4) \\
   & =\int_{\infty}^{(x_1,y_1)}\omega+\int_{\infty}^{(x_2,y_2)}\omega
     +\int_{\infty}^{(x_3,y_3)}\omega+\int_{\infty}^{(x_4,y_4)}\omega\
\end{aligned} \label{eq1.3}
  \\
\end{equation}
fills the whole space ${\mathbf C}^4$.  We denote the points in
${\mathbf C}^4$ by, for example, $u$ and $v$ and their natural
coordinates in ${\mathbf C}^4$ by the subscripts $(u_1,u_2,u_3,u_4)$,
$(v_1,v_2,v_3,v_4)$.  We denote the lattice generated by the integrals
of the basis (\ref{eq1.2}) along any closed paths on $C$ by $\Lambda$.
We denote the manifold ${\mathbf C}^4/\Lambda$, by $J$, the Jacobian
variety of $C$.  The projection from ${\mathbf C}^4$ to ${\mathbf
  C}^4/\Lambda$ is denoted by $\kappa$:
\begin{equation}
   \kappa:{\mathbf C}^4\rightarrow {\mathbf C}^4/\Lambda=J.
   \label{eq1.4}
\end{equation}
We have $\Lambda=\kappa^{-1}\big((0,0,0,0)\big)$.  We define for
$k=1$, $2$, $3$, $\dots$, the Abel map
\begin{equation}
\begin{aligned}
  \iota :\ \text{Sym}^k(C)&\rightarrow J \\
  (P_1,\cdots,P_k) &\mapsto \Big(\int_{\infty}^{P_1}\omega+\cdots+
  \int_{\infty}^{P_k}\omega\Big) \hbox{\rm mod}\,\Lambda,
   \end{aligned}\label{eq1.5}
\end{equation}
and denote its image by $W^{[k]}$.  Let
\begin{equation}
   [-1](u_1,u_2,u_3,u_4)=(-u_1,-u_2,-u_3,-u_4),
   \label{eq1.6}
\end{equation}
and 
\begin{equation}
   \Theta^{[k]}:=W^{[k]}\cup[-1]W^{[k]}.
   \label{eq1.7}
\end{equation} 
We call this  $\Theta^{[k]}$  the  $k$-th {\it standard theta subset}. 
In particular, if $k=1$, then (\ref{eq1.5}) gives an embedding of  $C$:
\begin{equation}
   \begin{aligned}
   \iota :\ &C\rightarrow J \\
   & P \mapsto \int_{\infty}^P\omega \ \hbox{\rm mod}\ \Lambda.
  \end{aligned}
   \label{eq1.8}
\end{equation}
We note that in contrast to the hyperelliptic case 
\begin{equation}
   \Theta^{[2]}\neq W^{[2]}, \quad \Theta^{[1]}\neq W^{[1]}.
\end{equation}  
On the embedded surface $\iota(C)=W^{[1]}$, we can take $u_4$ as a
local parameter at the origin $\iota(\infty)$.  Then we have (see
\cite{bg06,eel00}, for instance) Laurent expansions with respect to
$u_4$ as follows:
\begin{equation}
  u_1=\tfrac17{u_4}^7+\cdots,\quad u_2=\tfrac14{u_4}^4+\cdots,\quad 
  u_3=\tfrac12{u_4}^2+\cdots
  \label{u4expansion1}
\end{equation}
and also
\begin{equation}
   x(u)=\frac1{{u_4}^3}+\cdots, \quad y(u)=\frac1{{u_4}^5}+\cdots.
  \label{xyexpansion}
\end{equation}

As with the $(3,4)$ curve, we introduce a set of weights for the different
variables, as follows:
\begin{definition}
  We define a weight called the {\it Sato weight} for constants and
  variables appearing in our relations as follows.  The Sato weights
  of variables $u_1$, $u_2$, $u_3$, $u_4$ are $7$, $4$, $2$, $1$,
  respectively; the Sato weight of each coefficient $\lambda_j$ in
  (\ref{cyclic35}) is $15-3j$; while the Sato weights of $x(u)$ and
  $y(u)$ are $-5$ and $-3$, respectively.
  
  We note that the Sato weights of the variables $u_k$ are precisely
  the Weierstrass gap numbers of the Weierstrass gap sequence at
  $\infty$, whilst the Sato weights of $x(u)$ and $y(u)$ are
  Weierstrass non-gap numbers from the same sequence.
\end{definition}
All expressions in this paper are homogeneous with respect to this
weight.

\section{The sigma function} 
We  now construct the {\it sigma function} 
\begin{equation}
   \sigma(u)=\sigma(u_1,u_2,u_3, u_4)
   \label{sigma}
\end{equation} 
associated with  $C$  (see also \cite{bel97}, Chap.1). 
We choose a basis of cycles on $C$: 
\begin{equation}
   \alpha_i, \beta_j  \ (1\leqq i, j\leqq 4),
   \label{paths}
\end{equation} 
such that their intersection numbers are 
$$\alpha_i\cdot\alpha_j=\beta_i\cdot\beta_j=0,$$
$$\alpha_i\cdot\beta_j=\delta_{ij}.$$

Let  $Z$  and  $W$  be two indeterminates.  
We define
\begin{equation}
   \Omega\big((x,y),(z,w)\big)
   =\frac{1}{(x-z)\frac{\partial}{\partial y}f(x,y)}
    \sum_{k=1}^3y^{3-k}\bigg[\frac{f(Z,W)}{W^{3-k+1}}\bigg]_W\bigg|_{(Z,W)=(z,w)}
   \label{eq2.3} 
\end{equation} 
where  $[\ \ ]_W$  means removing the terms of negative powers 
with respect to  $W$. 

\begin{lemma} {\rm The fundamental 2-form of the second kind.}\ \
  \newline 
Let 
\begin{equation}
  \big((x,y),(z,w)\big)\mapsto R\big((z,w),(x,y)\big)\,\mathrm{d}z\,
  \mathrm{d}x
\end{equation}
be a $2$-form on $C\times C$ with only poles along the diagonal points
$\{((x,y),(x,y))\}\subset C\times C$, holomorphic elsewhere and
satisfying
\begin{equation}
  \lim_{x\to z}(x-z)^2\,R\big((z,w),(x,y)\big)=1. \label{polecond}
\end{equation}
For the differentials {\rm(\ref{eq1.2})}, for \,$\Omega$ above, and
for two variable points $(x,y)$ and $(z,w)$ on $C$, there exist second
kind differentials $\eta_j=\eta_j(x,y)$ $(j=1, \dots,4)$, having their
only pole at $\infty$, such that
\begin{equation}
  R\big((x,y),(z,w)\big)
  :=\frac{d}{dx}\Omega\big((x,y),(z,w)\big)
  +\sum_{j=1}^3\frac{\omega_j(x, y)}{dx}\frac{\eta_j(z, w)}{dz}, 
  \label{eq3.1.6} 
\end{equation} 
where the derivation\footnote{Since $x$ and $y$ are related, we do not
  use $\partial$.} is with respect to the variable point $(x,y)\in C$.
We further require that it satisfies the symmetry condition
\begin{equation}
  R\big((x,y),(z,w)\big)=R\big((z,w),(x,y)\big).
  \label{eq3.1.7} 
\end{equation}
Then the set of differentials $\{\eta_1$, $\eta_2$, $\eta_3$, $\eta_4 \}$
is determined uniquely modulo the space spanned by the $\omega_j$s of
{\rm(\ref{eq1.2})}.  The $2$-form obtained above is called 
{\rm(Klein's) fundamental $2$-form of the second kind}.
\end{lemma}
\begin{proof} 
The 2-form
\begin{equation}
\frac{d}{dx}\Omega\big((x,y),(z,w)\big)\mathrm{d}z\,\mathrm{d}x,\label{twoform}
\end{equation}
considered as a function of $(x, y)$, satisfies the condition on the poles;
indeed one can check that (\ref{twoform})
has only a second order pole at $(x, y) = (z,w)$ whenever $(z,w)$ is either
an ordinary point or a Weierstrass point; 
at infinity the expansion (\ref{xyexpansion}) should be used. 
However, the form (\ref{twoform}) has unwanted poles at infinity, when 
considered as a form in the (z,w)-variables. To restore the
symmetry required in (\ref{eq3.1.7}) we complement (\ref{twoform}) 
by the second term to obtain (\ref{eq3.1.6}) with
polynomials (x, y) which should be found from (2.15). That results in a 
system of linear equations for the $\eta_j(z,w)$, the coefficients of  
the $\omega(x,y)$, which is always solvable. As a result, the polynomials
$\eta_j(z,w)$ as well as $F((x, y), (z,w))$ are obtained explicitly, 
as in \cite{bg06}, pp.\ 3617--3618 (see also \cite{ba97}, around p.\ 194).
\end{proof}


It is easily seen that the $\eta_j$ above can be written as
\begin{equation}
 \eta_j(x,y)=\frac{h_j(x,y)}{\frac{\partial}{\partial y}f(x,y)}dx,
\end{equation}
where $h_j(x,y)\in\mathbf Q[\lambda_0,\cdots,
   \lambda_4][[x,y]]$, and $h_j$ is of homogeneous weight.
Straightforward calculations for the curve $C$ \eqref{cyclic35} then lead to
the following expressions:
\begin{eqnarray} \nonumber \eta_1(x,y) & = &y x \, ( 3 \lambda_3 + 7
  x^2 + 5x \lambda_4)\, \frac{\mathrm{d}x}{3x^2}, \\ \nonumber
  \eta_2(x,y) & =&2 y x \, \left( 2 x + \lambda_4
  \right)\,\frac{\mathrm{d}x}{3x^2}, \\ \nonumber 
\eta_3(x,y) & =&(2
  x^3 + x^2 \lambda_4 - \lambda_2)\, \frac{\mathrm{d}x}{3x^2}, \\
  \label{eq:trig 2nd kind} \eta_4(x,y) & =& xy\,
  \frac{\mathrm{d}x}{3x^2} .
\end{eqnarray}

We now define the period matrices by
\begin{equation}
   \left[\,\omega'  \ \omega''  \right]= 
\left[\int_{\alpha_i}\omega_j \ \ \int_{\beta_i}\omega_j\right]_{i,j=1,\dots,4},
   \,\,
   \left[\,\eta'  \ \eta''  \right]= 
\left[\int_{\alpha_i}\eta_j \ \ \int_{\beta_i}\eta_j\right]_{i,j=1,\dots,4}.
   \label{periods}
  \end{equation} 
We can combine these two matrices into 
\begin{equation}
   M=\left[\begin{array}{cc}\omega' & \omega'' \\ \eta' & \eta''
     \end{array}\right].
   \label{eq2.6}
\end{equation} 
The matrix  $M$ then satisfies 
\begin{equation}
   M\left[\begin{array}{cc} & -1_4 \\ 1_4 & \end{array}\right]{}^t {M}
   =2\pi\sqrt{-1}\left[\begin{array}{cc} & -1_4 \\ 1_4 &
     \end{array}\right].
   \label{Legendre}
\end{equation} 
This is the {\it generalized Legendre relation} (see (1.14) on p.\,11
of \cite{bel97}).  In particular, ${\omega'}^{-1}\omega''$ is a
symmetric matrix.  We know also that
\begin{equation}
   \text{Im}\,({\omega'}^{-1}\omega'') \ \ \ \text{is positive definite.}
   \label{positive_definiteness}
\end{equation}
Let
\begin{equation}
   \delta:=\left[\begin{array}{cc}\delta'\ \\
       \delta''\end{array}\right]\in \left(\tfrac12{\mathbf Z}\right)^{8}
   \label{eq2.9} 
\end{equation} 
be the theta characteristic which gives the Riemann constant with
respect to the base point $\infty$ and the period matrix $[\,\omega'\
\omega'']$ (\cite{mu85}, pp.163--166, \cite{bel97}, p.15, (1.18)).  By
looking at (\ref{eq1.2}), we see the canonical divisor class of $C$ is
given by $4\infty$.  Hence any theta characteristic is an element of
$\left(\tfrac12{\mathbf Z}\right)^{8}$ in this case.

We then define
\begin{equation}
\begin{aligned}
   \sigma(u)&=\sigma(u;M)=\sigma(u_1,u_2,u_3,u_4;M) \\
   &=c\,\text{exp}(-\tfrac{1}{2}u\eta'{\omega'}^{-1}\ ^t\negthinspace u)
   \vartheta\negthinspace
   \left[\delta\right]({\omega'}^{-1}\ ^t\negthinspace u;\ 
{\omega'}^{-1}\omega'') \\
   &=c\,\text{exp}(-\tfrac{1}{2}u\eta'{\omega'}^{-1}\ ^t\negthinspace u) \\
   &\hskip 10pt\times\sum_{n \in {\mathbf Z}^3} \exp \big[2\pi i\big\{
   \tfrac12 \ ^t\negthinspace (n+\delta'){\omega'}^{-1}\omega''(n+\delta')
   + \ ^t\negthinspace (n+\delta')(z+\delta'')\big\}\big], 
\end{aligned}
   \label{def_sigma}
\end{equation}
where $c$ is a constant depending only on the parameters of the curve,
$\{\lambda_0,\cdots,\lambda_4\}$, which we fix below. The series
(\ref{def_sigma}) converges for all $u\in{\mathbf C}^4$ because of
property (\ref{positive_definiteness}).

In what follows, for any given $u\in{\mathbf C}^4$, we denote by $u'$
and $u''$ the unique elements in ${\mathbf R}^4$ such that
\begin{equation}
   u=u'\omega'+u''\omega''.
   \label{eq2.12}
\end{equation}
Then for $u$, $v\in{\mathbf C}^4$, and $\ell$
($=\ell'\omega'+\ell''\omega''$) $\in\Lambda$, we define
\begin{align}
  L(u,v)    &:={}^t{u}(\eta'v'+\eta''v''),\nonumber \\
  \chi(\ell)&:=\exp[\pi\sqrt{-1}\big(2({}^t {\ell'}\delta''-{}^t
  {\ell''}\delta') +{}^t {\ell'}\ell''\big)] \ (\in \{1,\,-1\}).
   \label{eq2.13}
\end{align} 
The most important properties of $\sigma(u;M)$ can be expressed
as follows. 
\begin{lemma} 
\label{L2.14} 
For all $u\in{\mathbf C}^4$, $\ell\in\Lambda$,
and $\gamma\in\text{Sp}(8,{\mathbf Z})$, we have\,{\rm :} 
\begin{align}
\sigma(u+\ell;M) & =\chi(\ell)\sigma(u;M)\exp
L(u+\tfrac12\ell,\ell),\label{quasiperiodic}\\
\sigma(u;\gamma
M)&=\sigma(u;M),\label{modular}\\
  u\mapsto\sigma(u;M)& \text{ has zeroes of
order $1$ along } \Theta^{[3]}, \label{first-order}\\
\sigma(u;M)&=0 \iff u\in\Theta^{[3]}.\label{zeroes}
\end{align}
\end{lemma}

\begin{proof}
  The formula (\ref{quasiperiodic}) is a special case of the equation
  from \cite{ba98} (p.286, $\ell$.22).  The statement (\ref{modular})
  is easily shown by the definition of $\sigma(u)$ since $\gamma$
  corresponds to the choice of basis of cycles $\{\alpha_j,
  \beta_j\}_{j=1}^4$ (\ref{paths}), which are used to define the
  periods (\ref{periods}).  The statements (\ref{first-order}) and
  (\ref{zeroes}) are explained in \cite{ba98}, (p.252).  These facts
  are partially described also in \cite{bel97}, (p.12, Th. 1.1 and
  p.15).
\end{proof}
\begin{remark}%
\label{characterisation}
We fix a matrix $M$ satisfying (\ref{Legendre}) and
(\ref{positive_definiteness}).  The space of the solutions of
(\ref{quasiperiodic}) is a one dimensional space over ${\mathbf C}$,
because the Pfaffian of the Riemann form attached to $L(\, .\, )$ is 1
(see \cite{on98}, Lemma 3.1.2 and \cite{la82}, p.93, Th.3.1).  Hence,
such non-trivial solutions automatically satisfy
(\ref{modular})--(\ref{zeroes}).  In this sense, (\ref{quasiperiodic})
characterizes the function $\sigma(u)$ up to a constant factor. As a
corollary, since $\sigma(-u)$ also satisfies this condition, it
follows that $\sigma(u)$ must have definite odd or even parity.
\end{remark}

The constant $c$ of (\ref{def_sigma}) can be fixed as follows:
\begin{lemma}  
\label{L2.17}
The power series expansion of $\sigma(u)$ about $u=(0,0,0,0)$ with
respect to $u_1$, $u_2$, $u_3$ and $u_4$ has homogeneous Sato weight
\,$8$; its leading term is the Schur-Weierstrass polynomial $S(u)$
corresponding to the sequence of Sato weights of the $\{u_i\}$, which
is $\{7,4,2,1\}$ for any $(3,5)$ curve. Explicitly:
$$
S(u)=\left(\frac{1}{448}u_4^8 + u_2^2 + u_2 u_3 u_4^2 - \frac{1}{8}
  u_3^2 u_4^4 -\frac{1}{4} u_3^4 - u_1 u_4\right).
$$
The expansion is then of the form
\begin{equation*}
   \sigma(u)=\left(\frac{1}{448}u_4^8 + u_2^2 +  u_2 u_3 u_4^2 - 
   \frac{1}{8} u_3^2 u_4^4 -\frac{1}{4} u_3^4 - u_1 u_4\right)
   +(d^{\circ}({\lambda_0,\cdots,\lambda_4})\geqq 1).
\end{equation*}
In particular, by Remark \ref{characterisation}, all the higher 
monomials in this expansion must also be even, so here $\sigma(u)$ is 
an even function. 
\end{lemma}

\begin{proof}
  See the Section 2 of \cite{on05}. The essential part of this
  assertion is seen also by \cite{bel99}.
\end{proof}

In addition to the $\sigma$-function, we may also define $N$-th order
theta functions on ${\mathbf C}^4.$

\begin{definition}
  An $N$-th order theta function is any function $f(u)$ on ${\mathbf
    C}^4$ satisfying the same periodicity condition as $\sigma(u)^N$,
  that is:
\begin{align}
f(u+\ell) & =\chi(\ell)^N f(u)\exp(N L(u+\tfrac12\ell,\ell)).
\end{align}\label{Ntheta}
\end{definition}
Such functions will be used to construct Abelian functions below.

\section{Abelian functions} 
\begin{definition}
A meromorphic function $\mathfrak{P}(u)$ is called an Abelian function 
of $u\in \mathbb{C}^4$, with respect to the period lattice $\Lambda$ with
generators $\omega'$ and $\omega''$, if it is multiply periodic, that 
is, if 
\begin{equation} 
   \mathfrak{P}(u+\omega'n^T+\omega''m^T)=\mathfrak{P}(u) 
\end{equation}
for all integer vectors $n,m\in \mathbb{Z}$ wherever $\mathfrak{P}(u)$ exists.
\end{definition}
We note from  Definition \ref{Ntheta} that the quotient of two $N$-th order
theta functions must be Abelian.

To construct Abelian functions in terms of the $\sigma$-functions, we
first note that
$$
\sigma(u+v)\sigma(u-v)
$$ 
is a second order theta function in $u$. We then take derivatives with
respect to the parameter $v$, denoting:
\begin{equation}
  \Delta_i=\tfrac{\partial}{\partial v_i}
\end{equation}
for $u=(u_1,u_2,u_3,u_4)$ and $v=(v_1,v_2,v_3,v_4)$; then we define a
set of fundamental Abelian functions on $J$ by
\begin{equation}
 \wp_{ij}(u)=-\tfrac{1}{2 \sigma(u)^2}\Delta_i\Delta_j\,\sigma(u+v)
   \sigma(u-v)|_{v=0}
  =-\tfrac{\partial^2}{\partial u_i\partial u_j}\log\sigma(u).
\label{wpij}
\end{equation}
Evidently these functions are singular where $\sigma(u)=0$; this is on
the set $\Theta^{[3]}$.

For the benefit of the reader familiar only with the genus one case,
we should point out that the Weierstrass function $\wp(u)$ described
in eqn.\ (\ref{WP}) would be written as $\wp_{11}(u)$ in this notation.
Moreover, we define
\begin{equation}
   \wp_{ijk}(u)=\tfrac{\partial}{\partial u_k}\wp_{ij}(u),\ \
   \wp_{ijk\ell}(u)=\tfrac{\partial}{\partial u_{\ell}}\wp_{ijk}(u),
   \label{wp34}
\end{equation}
and so on for higher derivatives.  The functions (\ref{wpij}) and
(\ref{wp34}) are periodic functions because of (\ref{quasiperiodic}).
Moreover, following Baker \cite{ba07} and as generalised in
\cite{eemop06}, we define
\begin{equation}
  \begin{aligned}
    Q_{ijk\ell}(u)&=-\tfrac{1}{2\sigma(u)^2}\Delta_i\Delta_j\Delta_k\Delta_{\ell}
    \,\sigma(u+v)\sigma(u-v)|_{v=0}\\
    &= \wp_{ijk\ell}(u)-2(\wp_{ij}\wp_{k\ell}+\wp_{ik}\wp_{j\ell}+
    \wp_{i\ell}\wp_{jk})(u).
  \end{aligned}
\label{defQ}
\end{equation}

A short calculation shows that 
$$
Q_{ijk\ell} \in \varGamma(J,\mathcal{O}(2\Theta^{[3]})),
$$ 
that is the vector space of
meromorphic functions having at worst a double pole where $\sigma=0$,
but
$$
\wp_{ijk\ell} \in \varGamma(J,\mathcal{O}(4\Theta^{[3]})),
$$
having instead quadruple poles on the same set.  Indeed, we see
further that any expression of the Hirota form
\begin{equation}
  \begin{aligned}
    Q_{ijk\ell mn}(u)&=-\tfrac{1}{2\sigma(u)^2}\Delta_i\Delta_j
    \Delta_k\Delta_{\ell}\Delta_m\Delta_n
    \,\sigma(u+v)\sigma(u-v)|_{v=0},
  \end{aligned}
\label{defQ6}
\end{equation}
is an element of $\varGamma(J,\mathcal{O}(2\Theta^{[3]}))$, as is any
similar expression of higher even order in the $\Delta_i$; we note
that such expressions of odd order necessarily vanish. Further, these
may all be expressed as polynomials in the Kleinian $\wp_{ij}$ and
their derivatives.
 
Note that although the subscripts in $\wp_{ijk\ell}$ {\em do} denote
differentiation, the subscripts in $Q_{ijk\ell}$  do {\em not} denote
direct differentiation, and the latter notation is introduced for
convenience only.  This is important to bear in mind when we use
cross-differentiation, for example the $\wp_{ijk\ell}$ satisfy
\[
    \tfrac{\partial}{\partial u_m}\wp_{ijk\ell}(u) 
  = \tfrac{\partial}{\partial u_\ell}\wp_{ijkm}(u),
\]
whereas the $Q_{ijk\ell}$ do not.

The following formula involving the fundamental Kleinian
$\wp$-functions was derived in \cite{bel00} and was evaluated for the
case of the curve (\ref{eq1.2}) in \cite{bg06}. Here it is noted that
the unique Klein bidifferential may be written in two different ways,
either in terms of the second derivatives of $\ln(\sigma)$, or else in
terms of rational functions of the coordinates of points on the curve.

\begin{thm} \label{thm:trig PDE} For arbitrary $(x, y)$, and base
  point $\infty$ on $C$, an arbitrary set $S$ of $g=4$ distinct points
  $ \{ (x_1, y_1), \dots, (x_4, y_4) \} \in C^4$, and $(z,w)$ being
  any point of $S$, it follows that
\begin{equation}
  \sum_{i,j=1}^{4}{ \wp_{i,j} \left( \int_\infty^{t}{ \mathrm{d} {\mathbf u}}
      - \sum_{k=1}^{4}{ \int_\infty^{x_k}{ \mathrm{d} {\mathbf u}} } \right)} 
  \mathcal{U}_i(x,y)  \,  \mathcal{U}_j(z,w) = 
  \frac{F(x,y;z,w)}{(x-z)^2}, 
 \label{klein}
\end{equation}
where
\[
\mathcal{U}^{\mathrm{T}} (x,y) = (1, \, x, \, y, \, x^2) 
\]
the vector of numerators of the $\omega_i$, 
and $F$ is the symmetric function  
 \begin{align*} 
 F(&x,y;z,w) = 3 w^2 y^2 \\
 & +\left[ 2 z^3 x^2 + z^4 x + 3 \lambda_0 + \lambda_1 (2 z +x) 
   + \lambda_2 (z^2  +2 x z )+ \lambda_3 ( 3 z^2 x)  
   + \lambda_4 ( 2z^3 x +x^2 z^2 ) \right] y \\
 & + \left[ 2 x^3 z^2 + x^4 z + 3 \lambda_0 + \lambda_1 (2 x + z) 
   + \lambda_2 (x^2  +2 x z )+\lambda_3 ( 3 x^2 z) 
 + \lambda_4 ( 2x^3 z + x^2 z^2 ) \right] w
\end{align*}
which appears in the numerator of the second kind fundamental 2-form:
\[
\frac{F(x,y;z,w)}{(x-z)^2} 
\frac{\mathrm{d}x}{ f_y(x,y) } 
\frac{ \mathrm{d}z}{ f_{w}(z,w) }  = \mathrm{d} \Omega(x,y;z,w)
\]
\end{thm}
while $f(z,w)=0$ is the equation of the curve  $C$:
\[ 
f(z,w) = w^3 -( z^5+ \lambda_4 \, z^4 + \lambda_3 \, z^3 + \lambda_2
\, z^2 + \lambda_1 \, z + \lambda_0).
\]

Expanding (\ref{klein}) as one of the $x_i$ tends to infinity and
comparing the principal parts of the poles on both sides of the
relation, we find, in leading order, the solution of the Jacobi
inversion problem, first given explicitly for this curve in
\cite{bg06}, following \cite{eel00}.

\begin{thm}[Jacobi inversion formula]
  ~[\,\cite{bel97}, p. 32], \cite{bg06}, Let $C$ be the genus $4$ cyclic
  $(3,5)-$curve .  If $D=((x_1,y_1)+ (x_2,y_2) + (x_3,y_3) +
  (x_4,y_4))$ is a non-special divisor, $\mathrm{d} {\mathbf u}$ is
  the vector of holomorphic differentials, with period lattice
  $\Lambda$, then the Abel map is given by:
\[ 
u = \sum_{k=1}^{4} \int_{\infty}^{x_k}{ \mathrm{d}\mathbf{u}}
\qquad \mod \Lambda.
\] 

The Abel preimage of the point ${ u} \in {\mathbb C}^4 \,$ is then
given by the set
$$ 
\left\{ (x_1 , y_1), \dots, \, (x_4, y_4) \right\} \, \in \, 
\left( C \right) ^4,
$$
where $ \left\{ x_1, \dots,  x_4 \right\}$ are the zeros of the polynomial
\begin{eqnarray*} 
\mathcal{P}(x,y; { u} ) & =&2 x^4 +
 (\lambda_4 -  \wp_{444} -3  \wp_{34} ) x^3  \\ \nonumber
 && + 
\left( -  \wp_{34} \lambda_4 -  \wp_{44} \wp_{33} +  \wp_{444}  \wp_{34} +
 \wp_{34} \, \! ^2 -  \wp_{244} - \wp_{23} - \wp_{44}  \wp_{344} \right)
 x^2 \\ 
  && + 
\left( -\wp_{13} - \wp_{144} + \wp_{244} \wp_{34} + \wp_{23} \wp_{34} - \wp_{24}
\wp_{344} - \wp_{24} \wp_{33}   \right) x \\
&&+ \wp_{13} \wp_{34} -\wp_{14} \wp_{344} -\wp_{14} \wp_{33} + \wp_{144}
\wp_{34},
\end{eqnarray*}
and the coordinate $y_i$ of each point in $D$ is given by
\begin{equation} 
\wp_{14} + \wp_{24} x_i + \wp_{34} y_i + \wp_{44} x_i^2 - x_i y_i =0.
\end{equation}
\end{thm}

{\rm Expanding both sides of \eqref{klein} to higher order, we obtain
  successively higher order differential polynomials in the
  $\wp_{ij}$, which must vanish identically. These play an analogous
  role in the theory to Weierstrass' fundamental differential equation
  for the elliptic $\wp$, (\ref{WP}). One of these was given in
  \cite{bg06}, but more recent arguments based on the parity of
  $\sigma$ allow such equations to be greatly simplified.  One obtains
  sets of relations of third order, either linear or quadratic in the
  $\wp_{ijk}$, and others, of fourth order, linear in the
  $\wp_{ijkl}$.

The first few of the latter include:
\begin{align}
Q_{4444}& =-3\wp_{33}, \\
Q_{3444}& =3\wp_{24},\\
Q_{2344}& =-4\wp_{14}+ \wp_{22}-2 \lambda_{4} \wp_{24}.\\
Q_{1344}& =-\wp_{12}+2 \lambda_{4} \wp_{12}.
\end{align}
The simplest two quadratic relations are:
\begin{align}
\wp_{444}^2& = 4 \wp_{44}^3 - 4 \wp_{44}\wp_{33} +\wp_{34}^2
-4\wp_{23}+2 \lambda_4 \wp_{34} +\lambda_4^2 - 4 \lambda_3,\\
\wp_{344}^2& = 4 \wp_{34}^2 \wp_{44}^2 + 4 \wp_{24}\wp_{34}+\wp_{33}^2+4\wp_{14}.
\end{align}
Lists of all known relations of each class are given in section 
\ref{sec:relations on cyclic35}.}

\section{A basis of the space  $\varGamma(J, \mathcal{O}(2\Theta^{[2]}))$}

The vector space $\varGamma(J,\mathcal{O} (2\Theta^{[2]}))$ has
dimension $2^g=2^4=16$. All elements $\mathcal{V}(u)$ of this space may be
written as
\begin{equation}\label{E5.1}
\mathcal{V}(u) = \frac{f(u)}{\sigma(u)^2};
\end{equation}
here $\mathcal{V}(u)$ is an Abelian function and $f(u)$ is a second
order theta function.

We note that second-order theta functions with this period lattice form a
vector space, which by (\ref{E5.1}) is isomorphic to 
$\varGamma(J,\mathcal{O} (2\Theta^{[2]}))$. 
\begin{lemma}
We have the following basis
\begin{eqnarray*}
  \varGamma(J, \mathcal{O}(2\Theta^{[2]}))&
  = \mathbf{C} Q_{1144}
  \oplus{\mathbf C} \wp_{11}   
  \oplus{\mathbf C}\, Q_{1244}
  \oplus\mathbf{C}\, Q_{2233}\\
  &\oplus{\mathbf C} \wp_{12}
  \oplus{\mathbf C}\, Q_{1444}
  \oplus{\mathbf C} \wp_{13}
  \oplus{\mathbf C} \wp_{14}\\
  &\oplus{\mathbf C} \wp_{22}
  \oplus\mathbf{C}\, Q_{2444}
  \oplus{\mathbf C} \wp_{23}
  \oplus\mathbf{C} \wp_{24}\\
  &\oplus\mathbf{C} \wp_{33}
  \oplus\mathbf{C} \wp_{34}
  \oplus\mathbf{C} \wp_{44}
  \oplus\mathbf{C} 1  
   \label{basis}
 \end{eqnarray*}
\end{lemma}

\begin{proof} 
  We know the dimension of the space $\varGamma(J,\mathcal{O}
  (2\Theta^{[2]}))$ is $2^4=16$, by the Riemann-Roch theorem for
  Abelian varieties (see for example, \cite{mu85}, (pp.150--155),
  \cite{la82}, (p.99, Th.4.1).  Obviously, (\ref{wpij}, \ref{defQ})
  show that the functions on the right hand sides each belong to the
  space on the left hand side. It thus only remains to verify their
  linear independence, and it is sufficient to do this in the special
  case where all the $\lambda_j=0$, in which case $\sigma(u)$ reduces
  to the Schur-Weierstrass polynomial. We multiply each of the
  functions on the right hand side by $\sigma(u)^2$; this yields in
  each case a polynomial in $(u_1, u_2, u_3, u_4)$.  Then the
  functions on the right hand side are linearly independent; for all
  of these polynomials have different Sato weights, except for those
  from $\wp_{22}$ and $\wp_{14}$, and it is easy to check the
  independence of these two directly.
\end{proof}

\begin{cor}\label{L5.1}
A second-order theta function whose Taylor expansion at $(0,0,0,0)$ has a
leading term of Sato weight $16$, is a multiple of $\sigma(u)^2$.
\end{cor}
\begin{proof}
A second order theta function must be a linear combination of $\sigma(u)^2$
with other basis elements, all of Sato weight less than 16. These must be
absent if the leading order is 16.
\end{proof}
\begin{cor}\label{L5.2}
A second-order theta function whose Taylor expansion at $(0,0,0,0)$ has a
leading term of Sato weight greater than $16$, is identically zero.
\end{cor}
\begin{proof}
No non-trivial linear combination of these basis elements has leading weight
greater than 16.
\end{proof}

\section{Expansion of the $\sigma$-function} 

We need more terms of the power series expansion of $\sigma(u)$,
beyond those given in \cite{bg06}, in order to properly characterise
the curve. The terms of Sato weight 23 in $u_i$ are needed for this,
as these are the first containing any dependence on $\lambda_0$.
\begin{lemma} 
\label{L6.1} 
The function $\sigma(u)$ associated with {\rm (\ref{cyclic35})} has an
expansion of the following form\,{\rm :}
\begin{equation}
  \begin{aligned}
    \sigma(u_1,u_2,u_3,u_4)
    &=C_8(u_1,u_2,u_3,u_4)+C_{11}(u_1,u_2,u_3,u_4) \\
    &\hskip -20pt+C_{14}(u_1,u_2,u_3,u_4)+C_{17}(u_1,u_2,u_3,u_4)+\\
    &\hskip -20pt+C_{20}(u_1,u_2,u_3,u_4)+C_{23}(u_1,u_2,u_3,u_4)+\dots\\
    &\hskip -20pt+C_{5+3n}(u_1,u_2,u_3)+\dots
  \end{aligned}
\end{equation}
where each $C_{5+3n}$ is an even polynomial in the $u_i$, composed of
products of monomials in $u_i$ of total weight $5+3n$ multiplied by
monomials in $\lambda_i$ of total weight $-3n$.  The first few $C_k$
are given in \cite{bg06}.
\end{lemma}

\begin{proof} The proof is by construction (with heavy use of Maple)
  initially following the methods used in \cite{bg06}.  In particular,
  we require $C_8$ to be the Schur-Weierstrass polynomial, and then
  most of the higher terms are fixed by requiring the vanishing of the
  sigma function on the strata $\Theta^{[1]}$, $\Theta^{[2]}$ and
  $\Theta^{[3]}$. This determines $\sigma$ up to multiplication by an
  even analytic function equal to $1$ at the origin.  In addition we
  fix the remaining coefficients up to $C_{23}$ inclusive, by
  requiring that $\sigma$ satisfies (through the definitions
  (\ref{wpij}), (\ref{defQ})) the two equations:
  \begin{eqnarray}
  Q_{4444}& =&-3\wp_{33}, \\
  Q_{1344}& =&-\wp_{12}+2 \lambda_{4} \wp_{12}.
  \end{eqnarray}
  These equations, we recall, were found by expanding the Klein
  bidifferential and equating coefficients.
\end{proof}
\begin{remark}
  Other relations of this form, in which linear combinations of the
  $Q_{ijk\ell}$ are expressed as linear combinations of the $\wp_{ij}$
  could in principle be found by expanding the Klein bidifferential to
  (much) higher order, but it is much easier to derive them directly
  from the expansion of $\sigma$, as will be done below.
\end{remark}

\section{Equations satisfied by the Abelian functions associated with $C$.}
\label{sec:relations on cyclic35}  
We can use the $\sigma$-function expansion to identify various further
equations which the Abelian functions defined by (\ref{wpij}) and
(\ref{defQ}) must satisfy
\subsection{Four-index relations}
Such relations are the generalizations of
$\wp''=6\wp^2-\tfrac12g_2\wp$ in the cubic (genus 1) case.
\begin{proposition}
\label{L4index}
The 4-index functions $Q_{ijk\ell}$ associated with {\rm (\ref{cyclic35})}
satisfy the following relations\,
\begin{align*}
Q_{4444}& =-3\wp_{33}, \\
Q_{3444}& =3\wp_{24},\\
Q_{3344}& =-\wp_{23}+2\lambda_4 \wp_{34},\\
Q_{3334}& = -Q_{2444},\\
Q_{2344}& =-4\wp_{14}+ \wp_{22}-2 \lambda_{4} \wp_{24},\\
Q_{3333}& = 12\wp_{14}-3 \wp_{22},\\
Q_{2334}& =2 \wp_{13}+3 \lambda_3 \wp_{34},\\
Q_{1444}& =-\tfrac{1}{2}Q_{2333} +\tfrac{3}{2} \lambda_3 \wp_{33},\\
Q_{2244}& =-\tfrac{1}{3} Q_{2333}-\tfrac{2}{3}\lambda_4 Q_{3334}+2\lambda_3 
   \wp_{33},\\
Q_{1344}& = 2 \lambda_4 \wp_{14}-\wp_{12},\\
Q_{2234}& =-2\wp_{12} -2 \wp_{12}+4 \lambda_{4} \wp_{14}+3\lambda_3 \wp_{24}-2
\lambda_2 \wp_{44},\\
Q_{1334}& = -\tfrac{1}{2} Q_{2233} +2\lambda_4 \wp_{13} +\tfrac{3}{2} \lambda_3
\wp_{23} +2 \lambda_2 \wp_{34} + \lambda_4 \lambda_2,\\
Q_{1333}& = 3 Q_{1244}+\lambda_4 Q_{2333} -3 \lambda_4 \lambda_3 \wp_{33},\\
Q_{1234}& = -\wp_{11} +3 \lambda_3 \wp_{14} - \lambda_1 \wp_{44},\\
Q_{2223}& =6 \wp_{11} +6 \lambda_3 \wp_{14} +6 \lambda_3 \wp_{22} -6 \lambda_1
\wp_{44},\\
Q_{1223}& = -3\lambda_0+\lambda_4\lambda_1 +3\lambda_3 \wp_{13}+2\lambda_1
\wp_{34},\\
Q_{1144}& = -Q_{1224}-\tfrac{1}{2}\lambda_3 Q_{2333} +\tfrac{3}{2}\lambda_3^2
\wp_{33}+ 3 \lambda_1 \wp_{33},\\
Q_{1134}& = \tfrac{2}{3}\lambda_4 Q_{2223}+(2 \lambda_2 -4\lambda_3 \lambda_4)
\wp_{14}\\
& \qquad -\lambda_1 \wp_{24}+4 \lambda_1 \lambda_4 \wp_{44} -2\lambda_3 \lambda_4
\wp_{22}+4 \lambda_4 \wp_{11},\\
Q_{1223}& = -2\lambda_4 \wp_{11} +3\lambda_3\wp_{12}+4\lambda_2\wp_{14}
 -2\lambda_1\wp_{24} -6\lambda_0\wp_{44},\\
 Q_{1222}& = 6( \lambda_0+\lambda_4 \lambda_1) \wp_{33}
    -\lambda_3 Q_{1333}+ 2 \lambda_1 Q_{2444},\\
\ldots.&&
\end{align*}
\end{proposition}
\begin{proof}
  Since any 4-index function $Q_{ijk\ell}$ belongs to
  $\varGamma(J,\mathcal{O}(2\Theta^{[2]}))$, it must be in the span of the
  basis (\ref{basis}) - hence a relation must hold of the form:
\[
Q_{ijk\ell} = \text{ linear function of various } \wp_{mn} \text{ and
   other } Q_{opqr},
\]
with $\lambda$-dependent coefficients. Thus, it is just a matter of
enumerating all the possible terms on the right-hand side which have
the same weight as $Q_{ijk\ell}$.  Using the method of undetermined
coefficients, we insert the $\sigma$ expansion truncated to an
appropriate weight in the $\lambda_i$ and solve for the unknown
(constant) coefficients.  This requires the use of Maple but is quite
efficient as the pole on either side is only of order 2 in  $\sigma$.
\end{proof}
\begin{remark} For completeness we give the explicit formulae for the
  $Q$ in terms of the $\wp$
\begin{align*}
&  Q_{ijkk} = \wp_{ijkk} - 2 \wp_{ij}\wp_{kk}-4\wp_{ik}\wp_{jk},\ \ 
&& Q_{iikk} = \wp_{iikk} - 2 \wp_{ii}\wp_{kk}-4\wp_{ik}^2,\\
&  Q_{ikkk} = \wp_{ikkk} - 6\wp_{ik}\wp_{kk},\ \
&& Q_{kkkk} = \wp_{kkkk} - 6\wp_{kk}^2.
\end{align*}
\end{remark}
\begin{remark}     
  The complete set of such relations for the cyclic $(3,4)$ curve was
  given in \cite{eemop06}. As far as we know, the above incomplete set
  of relations for the cyclic $(3,5)$ curve is new, and a comparison
   is of interest.
\end{remark}

\begin{remark}
The first equation (\ref{L4index}) above is:
\[
\wp_{4444}=6\,  \wp_{44} ^{2}-3\,\wp_{33},
\] 
This relation, after differentiating twice with respect to $u_4$,
becomes the Boussinesq equation for $\wp_{44}$; $u_4$ and $u_3$
respectively play the roles of the space and time variables here. The
connection between the Boussinesq equation and cyclic trigonal curves
is well-established (see \cite{bel00,eel00}).
\end{remark}

\subsection{Linear three-index relations}
\begin{proposition}
\label{L3index}
The $3$-index functions $\wp_{ijk}$ associated with {\rm
  (\ref{cyclic35})} satisfy a number of relations linear in these
functions.  These have no analogue in the genus 1 case.  For example
in decreasing weight, starting at -6, we have
\begin{align*}
\wp_{333} & =  2 \wp_{44} \wp_{344} -2 \wp_{34} \wp_{444} -\wp_{244},\\
\wp_{234} & = \tfrac12\wp_{34} \wp_{344} - \wp_{334}  \wp_{44}
  +\tfrac12 \wp_{33} \wp_{444} +\tfrac12 \lambda_4 \wp_{344},\\
\wp_{233} & =-\wp_{33} \wp_{344}  -\tfrac32 \wp_{444} \wp_{24}
  +\tfrac12 \wp_{334} \wp_{34}+\tfrac32 \wp_{244} \wp_{44}+\tfrac12 
\lambda_4\wp_{334}
  +\tfrac12 \wp_{333} \wp_{44},\\
\wp_{144} &  = -\tfrac12 \wp_{334} \wp_{33} + \tfrac12 \wp_{333} \wp_{34}
  +\wp_{344}\wp_{24}-\tfrac12 \wp_{34}\wp_{244},\\
\wp_{134} & = \wp_{234} \wp_{34} -\wp_{24}\,\wp_{334} 
  + \tfrac12\,\wp_{33}\wp_{244}-\tfrac12 \,\wp_{344}\,\wp_{23},\\
\wp_{133} & = \tfrac12 \wp_{333} \wp_{24}-\wp_{33} \wp_{234}
  -\tfrac12 \wp_{23}\wp_{334}+\wp_{34}\wp_{233}-3\wp_{444}\wp_{14}+
3\wp_{144}\wp_{44},\\
\wp_{124} & = -\wp_{134} \wp_{44}- \tfrac12 \wp_{144} \wp_{34}+ \wp_{14}
  \wp_{344} +\tfrac12  \wp_{13} \wp_{444}+\tfrac12 \lambda_4 \wp_{144},\\
\wp_{134} \wp_{34} & = -\tfrac12 \wp_{33} \wp_{144}
  +\tfrac12 \wp_{344}\wp_{13} +  \wp_{334} \wp_{14},\\
  \wp_{114} &= -\tfrac12 \wp_{144}\wp_{23} -\wp_{134} \wp_{24} +
  \wp_{234} \wp_{14}+\tfrac12 \wp_{244}\wp_{13},\\
  \wp_{111} & =\tfrac23 \wp_{22} \wp_{123}+\tfrac13 \wp_{23} \wp_{122}
  +\lambda_3 \wp_{114}-\lambda_1 \wp_{144} - \tfrac13 \wp_{13} \wp_{222}\\
 & \qquad-\tfrac23 \wp_{223}\wp_{12} - \tfrac23 \lambda_2 \wp_{124}
  +\tfrac13 \lambda_1 \wp_{224} +\lambda_0 \wp_{244} +\tfrac13 \lambda_4
  \wp_{112},
\end{align*}
\end{proposition}
\begin{proof}
  These can be calculated directly by expressing the equations in
  Proposition \ref{L4index} in terms of $\wp_{ijk\ell}$ and $\wp_{mn}$
  functions, then using cross differentiation on suitably chosen pairs
  of equations.  For example, the first relation above for $\wp_{333}$
  comes from
\[
\frac{\partial}{\partial u_3}\wp_{4444} - \frac{\partial}{\partial
   u_4}\wp_{3444}=0.
\]
In principle the relations could be checked by inserting the $\sigma$
expansion to the required level, but this is method is much slower and more
cumbersome since the terms involve poles of order 3.
\end{proof}
\subsection{Quadratic three-index relations}
\begin{proposition}
\label{L33index}
Quadratic expressions in the $3$-index functions $\wp_{ijk}$
associated with {\rm (\ref{cyclic35})} can be expressed in terms of {\rm
  (}at most cubic{\rm )} relations in the $\wp_{mn}$ and
$\wp_{ijk\ell}$.  For example we have the following
\begin{align*}
   \wp_{444}^2 & = 4 \wp_{44}^3-4 \wp_{44} \wp_{33}+\wp_{34}^2-4
   \wp_{23}
   +2 \lambda_{4} \wp_{34}+\lambda_{4}^2-4 \lambda_{3},  \\
   \wp_{344} \wp_{444} & = 4 \wp_{34} \wp_{44}^2+6 \wp_{24} \wp_{44}
   -\wp_{33} \wp_{34}- \wp_{33} \lambda_{4}-\tfrac23 \wp_{2444} , \\
   \wp_{344}^2 & = 4 \wp_{34}^2 \wp_{44} + 4 \wp_{24} \wp_{34} +
   \wp_{33}^2 +4 \wp_{14} , \\
   \wp_{334} \wp_{444} & = 2 \wp_{34}^2 \wp_{44}-\wp_{24} \wp_{34}-
   2 \wp_{33}^2-4 \wp_{14}- 2 \wp_{44} \wp_{23}+ 2 \wp_{22}\nonumber \\
   & \qquad-\wp_{24} \lambda_{4}+
   2 \wp_{44} \lambda_{4} \wp_{34}+ 2 \wp_{33} \wp_{44}^2,\\
   \wp_{334} \wp_{344} & = 2 \wp_{44} \wp_{33} \wp_{34}+2 \wp_{34}^3
   -2 \wp_{23} \wp_{34}+\wp_{24} \wp_{33}-2 \wp_{13}+2 \lambda_4
   \wp_{34}^2,
\\
   \wp_{333} \wp_{444} & = 6 \wp_{44} \wp_{33} \wp_{34}-2 \wp_{34}^3
   +7 \wp_{34} \wp_{23}+2 \wp_{33} \wp_{24}+2 \wp_{13}-4 \lambda_4
   \wp_{34}^2
   \nonumber\\
   & \qquad +2 \lambda_2+6 \wp_{34} \lambda_3+\lambda_4 \wp_{23}-2
   \lambda_4^2 \wp_{34}
   -2 \lambda_4 \wp_{44} \wp_{33}+12 \wp_{24} \wp_{44}^2\nonumber \\
   & \qquad -2 \wp_{44} \wp_{2444}, \\
   \wp_{444} \wp_{244} & = \tfrac23 \wp_{44} \wp_{2444}-2 \wp_{13}
   +\wp_{34} \wp_{23}-2 \wp_{33} \wp_{24}-2 \lambda_2+2 \wp_{34}
   \lambda_3\nonumber\\
& \qquad -\lambda_4 \wp_{23},
 \\   
   \wp_{334}^2 & = 4 \wp_{33} \wp_{34}^2+8 \wp_{34} \wp_{24} \wp_{44}
   -\tfrac43 \wp_{2444} \wp_{34}+\wp_{24}^2-\tfrac43 \wp_{1444}\nonumber\\
& \qquad
   +4 \wp_{44} \wp_{14}, \\
   \wp_{344} \wp_{333} & = 2 \wp_{33} \wp_{34}^2-\wp_{33} \wp_{23} +2
   \wp_{33} \lambda_4 \wp_{34}+2 \wp_{33}^2 \wp_{44}-4 \wp_{34}
   \wp_{24} \wp_{44}\nonumber\\
& \qquad  +\tfrac23 \wp_{2444} \wp_{34}-2
   \wp_{24}^2+\tfrac23 \wp_{1444}+4 \wp_{44} \wp_{14},
   \\
   \wp_{344} \wp_{244} & = \wp_{33} \wp_{23}+\tfrac23 \wp_{2444}
   \wp_{34} +2 \wp_{24}^2-\tfrac23 \wp_{1444}+4 \wp_{44} \wp_{14},\\
   \wp_{444} \wp_{234} & = 4 \wp_{34} \wp_{24} \wp_{44} -\tfrac13
   \wp_{2444} \wp_{34}-2 \wp_{33} \wp_{23}+4 \lambda_4 \wp_{24}
   \wp_{44} -\tfrac13 \lambda_4 \wp_{2444}\nonumber\\
& \qquad -2 \wp_{33} \lambda_3+4
   \wp_{44} \wp_{14} +2 \wp_{23} \wp_{44}^2-2 \wp_{2,2} \wp_{44},\\
   \wp_{344} \wp_{234} & = 2 \wp_{14} \lambda_4 +2 \wp_{44} \wp_{33}
   \wp_{24}+2 \wp_{34} \wp_{23} \wp_{44}+2 \wp_{44} \wp_{13} +2
   \wp_{14} \wp_{34}\nonumber\\
& \qquad +2 \wp_{24} \wp_{34}^2+2 \wp_{34} \lambda_4
   \wp_{24} -\tfrac13 \wp_{33} \wp_{2444},\\
   \wp_{444} \wp_{233} & = 6 \wp_{44} \wp_{34} \lambda_3 -2 \wp_{44}
   \lambda_2+4 \wp_{34} \wp_{23} \wp_{44}-4 \wp_{34} \lambda_4
   \wp_{24} -2 \wp_{24} \wp_{34}^2\nonumber\\
& \qquad -2 \wp_{14} \wp_{34}+\wp_{34}
   \wp_{2,2}-2 \wp_{44} \wp_{13} -2 \wp_{44} \wp_{33}
   \wp_{24}+\lambda_4 \wp_{2,2}\nonumber \\
& \qquad -2 \lambda_4 \wp_{23} \wp_{44}
   +\tfrac23 \wp_{33} \wp_{2444}+6 \wp_{24} \wp_{23}-2 \wp_{24}
   \lambda_4^2 +6 \wp_{24} \lambda_3\nonumber \\
& \qquad-2 \wp_{14} \lambda_4,
   \\
   \wp_{334} \wp_{244} & = 2 \wp_{34} \lambda_4 \wp_{24} +2 \wp_{24}
   \wp_{34}^2+2 \wp_{14} \wp_{34}-2 \wp_{34} \wp_{22}+2 \wp_{12} -2
   \wp_{44} \wp_{13}\nonumber\\
& \qquad -2 \wp_{44} \wp_{33} \wp_{24}+\tfrac23 \wp_{33}
   \wp_{2444} -\wp_{24} \wp_{23}-2 \wp_{14} \lambda_4,
   \\
   \wp_{334} \wp_{333} & = 6 \wp_{14} \wp_{34}-2 \wp_{34} \wp_{22} -2
   \wp_{12}-2 \wp_{44} \wp_{13}+4 \wp_{44} \wp_{33} \wp_{24}\nonumber\\
& \qquad 
   -\tfrac23 \wp_{33} \wp_{2444}+\wp_{24} \wp_{23}+2 \wp_{14}
   \lambda_4 +4 \wp_{34} \wp_{33}^2, \\
   \wp_{333}^2 & = 2 \wp_{2233}-4 \lambda_4 \lambda_2 +4 \lambda_1 -8
   \wp_{33} \wp_{22} -4 \wp_{34} \wp_{13}-6 \wp_{23} \lambda_3 \nonumber\\
& \qquad-4
   \wp_{34} \lambda_2 -4 \lambda_4 \wp_{13} -7 \wp_{23}^2 +16 \wp_{14}
   \wp_{33} +4 \wp_{33}^3,
   \\
   \wp_{444} \wp_{144} & = \tfrac23 \wp_{33} \wp_{22}+\wp_{34}
   \wp_{13} +\tfrac23 \wp_{44} \wp_{1444}+\tfrac43 \wp_{23}^2-2
   \wp_{33} \wp_{14}-\tfrac13 \wp_{2233}\nonumber\\
& \qquad +\tfrac13 \lambda_4
   \wp_{13}+\wp_{23} \lambda_3+\tfrac43 \wp_{34} \lambda_2 -\tfrac43
   \lambda_1+\tfrac23 \lambda_4 \lambda_2, \\
   \wp_{244}^2 & = -4 \wp_{44} \wp_{24}^2-\tfrac43 \wp_{33} \wp_{22}
   -\tfrac53 \wp_{23}^2+\tfrac43 \wp_{24} \wp_{2444}+\tfrac23
   \wp_{2233} -\tfrac83 \lambda_4 \wp_{13}\nonumber\\
& \qquad-2 \wp_{23} \lambda_3
   +\tfrac43 \wp_{34} \lambda_2-\tfrac43 \lambda_1-\tfrac43 \lambda_4
   \lambda_2,
   \\
   \wp_{244} \wp_{333} & = -\tfrac43 \wp_{24} \wp_{2444} +(\tfrac83
   \wp_{13}+2 \wp_{34} \wp_{23}) \lambda_4-\tfrac23 \wp_{2233}
   +\tfrac53 \wp_{23}^2+\tfrac43 \wp_{33} \wp_{22}\nonumber\\
& \qquad+4 \wp_{34} \wp_{13}
   +4 \wp_{34} \wp_{33} \wp_{24}-2 \wp_{34}^2 \wp_{23}-\tfrac43
   \wp_{44} \wp_{1444} +8 \wp_{44} \wp_{24}^2\nonumber\\
& \qquad+2 \wp_{44} \wp_{33}
   \wp_{23}+8 \wp_{14} \wp_{44}^2 +2 \wp_{23} \lambda_3+\tfrac83
   \wp_{34} \lambda_2+\tfrac43 \lambda_1 +\tfrac43 \lambda_4
   \lambda_2\nonumber\\
& \qquad-4 \wp_{34}^2 \lambda_3,\\
   \wp_{233} \wp_{344} & = \tfrac23 \wp_{24} \wp_{2444} +\tfrac13
   \wp_{2233}-\tfrac43 \wp_{23}^2-\tfrac53 \wp_{33} \wp_{22}-\wp_{23}
   \lambda_3 +\tfrac23 \wp_{34} \lambda_2\nonumber\\
& \qquad+2 \wp_{34}^2
   \wp_{23}-\tfrac43 \wp_{44} \wp_{1444} -4 \wp_{44} \wp_{24}^2+2
   \wp_{44} \wp_{33} \wp_{23}+2 \wp_{34}^2 \lambda_3 \nonumber\\
& \qquad+(2 \wp_{33} \wp_{24}-\tfrac43
   \wp_{13}) \lambda_4+8 \wp_{14} \wp_{44}^2
   +\tfrac43 \lambda_1-\tfrac23 \lambda_4 \lambda_2,
   \\
   \wp_{234} \wp_{334} & = -\tfrac13 \wp_{24} \wp_{2444}+\tfrac13
   \wp_{2233} -\tfrac43 \wp_{23}^2-\tfrac23 \wp_{33} \wp_{22}+2
   \wp_{33} \wp_{14}-\wp_{23} \lambda_3\nonumber\\
& \qquad +\tfrac23 \wp_{34}
   \lambda_2-\tfrac43 \lambda_4 \wp_{13}+2 \wp_{34} \wp_{33} \wp_{24}
   +2 \wp_{34}^2 \wp_{23}+\tfrac23 \wp_{44} \wp_{1444}\nonumber\\
& \qquad+2 \wp_{44}
   \wp_{24}^2 +2 \wp_{34}^2 \lambda_3-4 \wp_{14} \wp_{44}^2+\tfrac43
   \lambda_1 -\tfrac23 \lambda_4 \lambda_2,
   \\
   \wp_{224} \wp_{444} & = -\tfrac43 \lambda_1 -\tfrac13 \lambda_4
   \lambda_2+2 \wp_{34} \lambda_4 \wp_{23} +\lambda_4 \wp_{34}
   \lambda_3+2 \wp_{44} \wp_{33} \lambda_3 \nonumber\\
& \qquad-2 \lambda_4 \wp_{33}
   \wp_{24}+2 \wp_{22} \wp_{44}^2 +\tfrac23 \wp_{44} \lambda_4
   \wp_{2444}-4 \lambda_4 \wp_{24} \wp_{44}^2 -\tfrac23 \wp_{23}^2\nonumber\\
& \qquad-2
   \wp_{33} \wp_{14}+2 \wp_{34} \wp_{13}+\tfrac23 \wp_{44} \wp_{1444}
   -\wp_{34}^2 \lambda_3-4 \wp_{14} \wp_{44}^2\nonumber\\
& \qquad-\tfrac13 \wp_{2233} +6
   \wp_{44} \wp_{24}^2+\tfrac23 \wp_{33} \wp_{22}-\tfrac23 \wp_{24}
   \wp_{2444} -\tfrac23 \lambda_4 \wp_{13}\nonumber\\
& \qquad-\wp_{23} \lambda_3+\tfrac73
   \wp_{34} \lambda_2,
\\\wp_{144}\wp_{344} &= \tfrac23 \wp_{1444}\wp_{34}+\wp_{13}\wp_{33}+
2\wp_{14}\wp_{24}
\end{align*}
\end{proposition}

\begin{proof}
  Such relations can be found in three different ways.  One is to
  multiply one of the linear three-index $\wp_{ijk}$ relations above
  by another $\wp_{ijk}$ and substitute for previously calculated
  $\wp_{ijk}\wp_{\ell mn}$ relations of higher weight.  A second
  approach is to take a derivative of one of the linear three-index
  $\wp_{ijk}$ relations above and to substitute the known linear
  four-index $\wp_{ijk\ell}$ and previously calculated
  $\wp_{ijk}\wp_{\ell mn}$ relations.  The third method is to use the
  method of undetermined coefficients: write down all possible terms
  which are cubic or less in the basis functions (\ref{basis}) of the
  correct weight.  This last method is guaranteed to work but is much
  more time-consuming than the first two methods (when they work).  At
  weight -15 this method requires some use of Distributed Maple on a
  multi-processor cluster.
\end{proof}
\begin{remark}
  These relations are the generalizations of the familiar result $
  \left(\wp'\right)^2 = 4 \wp^3-g_2\wp-g_3$ in the genus 1 theory.
\end{remark}

\begin{remark}
  It is evident that the left hand sides of the above relations are
  not independent.  For instance from the first three of them, we may
  find two different expressions for $\wp_{444}^2\,\wp_{344}^2$.
  Identifying these we obtain a quadratic four-index relation:
$$( 4 \wp_{44}^3-4 \wp_{44} \wp_{33}+\wp_{34}^2-4 \wp_{23}
   +2 \lambda_{4} \wp_{34}+\lambda_{4}^2-4 \lambda_{3})
   (4 \wp_{34}^2 \wp_{44} + 4 \wp_{24} \wp_{34} + \wp_{33}^2 +4 \wp_{14})$$
$$-(4 \wp_{34} \wp_{44}^2+6 \wp_{24} \wp_{44}
   -\wp_{33} \wp_{34}- \wp_{33} \lambda_{4}-\tfrac23 \wp_{2444} 
   )^2 =0.$$
We note that the linear 4-index relations above do not include one involving
$\wp_{2444}$ alone. 
\end{remark}

\section{The two-term addition theorem} 

\begin{theorem} 
\label{Two-term}
The sigma function associated with {\rm (\ref{cyclic35})} satisfies the
following two-term addition formula\,{\rm :}
\begin{equation}
  \begin{aligned}
   &\frac{\sigma(u+v)\sigma(u-v)}{\sigma(u)^2\sigma(v)^2}
    = -\wp_{11}(u)\wp_{44}(v) + \wp_{12}(u)\wp_{24}(v)-
    \tfrac34 \wp_{14}(u)\wp_{22}(v)\\
   &\quad +\tfrac13 \wp_{13}(u)Q_{2444}(v) +\tfrac{1}{12} \wp_{14}(u)Q_{3333}(v)
     +\tfrac16 \wp_{23}(u) Q_{2333}(v) + \tfrac13 \wp_{33}(u) Q_{1334}(v)\\
   &\quad -\tfrac13 \wp_{34}(u) Q_{1333}(v)-\tfrac{1}{12} Q_{2222}(u)
     - \tfrac13 \lambda_4 Q_{1333}(u) + \tfrac16 \lambda_3 Q_{2333}(u)\\
   & \quad
     -\tfrac12 \lambda_3 \wp_{23}(u)\wp_{33}(v) +\tfrac13 \lambda_2 Q_{2444}(u) +\left(\tfrac13 \lambda_1 
     +\lambda_4 \lambda_2 -\tfrac34 \lambda_3^2\right)\wp_{33}(u)
     + (u\Leftrightarrow v).\label{2Term}
  \end{aligned}
\end{equation}
\end{theorem}

\begin{proof}
  Firstly, we notice that the left hand side is an even function with
  respect to $(u,v)\mapsto ([-1]u,[-1]v)$, and is a symmetric Abelian function
  in $u$ and   $v$.  We note that it has poles of order 2 along
  $(\Theta^{[2]}\times J)\cup(J\times\Theta^{[2]})$ but nowhere else.
  Moreover it is of Sato weight $-16$. It is thus expressible as a
  symmetric bilinear combination of the basis functions in
  (\ref{basis}), with each pair of functions having total weight
  $-16+3n$, where the coefficients are either absolute constants, (if
  $n=0$) or else are homogeneous polynomials (of weight $-3n$) in the
  $\lambda_i$. These undetermined coefficients were then found by
  substituting in the expansion of $\sigma$, and equating
  coefficients.
\end{proof}

\begin{remark} 
\label{R7.4}
By applying 
\begin{equation}
   \tfrac12\frac{\partial}{\partial u_i}
   \Big(
   \frac{\partial}{\partial u_j}
  +\frac{\partial}{\partial v_j}
   \Big)\log
\end{equation}
to \ref{2Term}, we have $\wp_{ij}(u+v)- \wp_{ij}(u)$ on the left hand
side, and have a rational expression of several
$\wp_{ij\cdots\ell}(u)$s and $\wp_{ij\cdots\ell}(v)$s on the right
hand side.  Hence, we have algebraic addition formulae for
$\wp_{ij}(u)$s.
\end{remark}
\begin{remark}
  We note that the left-hand side of the addition formula has a zero
  wherever $u=v$. By putting $v=u-w$ and letting $w \rightarrow
  (0,0,0,0)$, the left hand side has leading term
$$
\frac{\sigma(u+v)}{\sigma(u)^2\sigma(v)^2}(w_2^2-w_1 w_4) +O(w^3).
$$
On expanding the left and right hand sides, we can match powers of
$w$, getting at zeroth order:
\begin{cor}
\begin{equation}
  \begin{aligned}
    &- \wp_{11}(u)\wp_{44}(u) + \wp_{12}(u)\wp_{24}(u)-
    \tfrac34 \wp_{14}(u)\wp_{22}(u) +\tfrac13 \wp_{13}(u)Q_{2444}(u)\\
    &\quad +\tfrac{1}{12} \wp_{14}(u)Q_{3333}(u)+\tfrac16 \wp_{23}(u)
    Q_{2333}(u) + \tfrac13 \wp_{33}(u) Q_{1334}(u)\\
    &\quad -\tfrac13 \wp_{34}(u) Q_{1333}(u)-\tfrac{1}{12}
    Q_{2222}(u)- \tfrac13 \lambda_4 Q_{1333}(u) + 
      \tfrac16 \lambda_3 Q_{2333}(u)\\
    &\quad-\tfrac12 \lambda_3 \wp_{23}(u)\wp_{33}(u)+\tfrac13
    \lambda_2 Q_{2444}(u)+\left(\tfrac13 \lambda_1 +\lambda_4
      \lambda_2 - \tfrac34 \lambda_3^2\right)\wp_{33}(u) =0.
  \end{aligned}
\end{equation}
\end{cor}
At first order, we obtain four equations, for instance:
\begin{cor}
\begin{equation}
  \begin{aligned}
    &- \wp_{11}(u)\wp_{144}(u) + \wp_{12}(u)\wp_{124}(u)-\tfrac34
    \wp_{14}(u)\wp_{122}(u)\\
    & \quad+\tfrac13 \wp_{13}(u)(Q_{2444}(u))_1 +\tfrac{1}{12} \wp_{14}(u)
    (Q_{3333}(u))_1+\tfrac16 \wp_{23}(u) (Q_{2333}(u))_1 \\
    &\quad+ \tfrac13 \wp_{33}(u) (Q_{1334}(u))_1-\tfrac13 \wp_{34}(u)
    (Q_{1333}(v))_1 + \tfrac16 \lambda_3 Q_{2333}(u)\\
    &\quad-\tfrac12 \lambda_3 \wp_{23}(u)\wp_{133}(u) -
    \wp_{111}(u)\wp_{44}(u)
    + \wp_{112}(u)\wp_{24}(u)-\tfrac34 \wp_{114}(u)\wp_{22}(u) \\
    &\quad+\tfrac13 \wp_{113}(u)Q_{2444}(u)
    +\tfrac{1}{12} \wp_{114}(u)Q_{3333}(u)+\tfrac16 \wp_{123}(u) Q_{2333}(u)\\
    &\quad + \tfrac13 \wp_{133}(u) Q_{1334}(u)
    -\tfrac13 \wp_{134}(u) Q_{1333}(u)-\tfrac{1}{12} (Q_{2222}(u))_1\\
    &\quad - \tfrac13 \lambda_4 (Q_{1333}(u))_1 + \tfrac16 \lambda_3
    (Q_{2333}(u))_1-\tfrac12 \lambda_3
    \wp_{123}(u)\wp_{33}(u)\\
    &\quad +\tfrac13 \lambda_2 (Q_{2444}(u))_1+\left(\tfrac13 \lambda_1
      +\lambda_4 \lambda_2 -\tfrac34 \lambda_3^2\right)\wp_{133}(u)=0.
  \end{aligned}
\end{equation}
\end{cor}
At second order, we obtain several identities, in particular a
``double-angle'' sigma formula:
\begin{cor}
\begin{equation}
\begin{aligned}
  \nonumber \frac{\sigma(2u)}{\sigma(u)^4}& = \tfrac16
  \wp_{33}\wp_{122334}-\tfrac{7}{12}\wp_{1224}\wp_{33}^2-\tfrac34
  \lambda_3  \wp_{23}\wp_{2233}+\tfrac12 \wp_{222}^2-\tfrac{1}{24}\wp_{222222}\\
  &\quad+\tfrac16 \lambda_2 \wp_{222444}+\tfrac12 \lambda_1
  \wp_{2233}-\tfrac38 \lambda_3^2 \wp_{2233}+\tfrac16
  \wp_{1223}\wp_{2444}+\tfrac16 \wp_{2233}\wp_{1334}
  \\
  &\quad+\tfrac{1}{24}\wp_{1224}\wp_{3333}+\tfrac{1}{24}
  \wp_{14}\wp_{223333}+\tfrac16
  \wp_{13}\wp_{222444}-\wp_{23}\wp_{233}\wp_{223}\\
  &\quad+\tfrac{1}{12}\wp_{23}\wp_{222333}-\tfrac12 \wp_{2233}\wp_{23}^2
  +\tfrac13
  \wp_{33}\wp_{2234}\wp_{13}+\tfrac{1}{12}\wp_{2223}\wp_{2333}
  -\tfrac12\wp_{14}\wp_{233}^2\\
  &\quad-\lambda_2 \wp_{44}\wp_{2224} -2 \lambda_2 \wp_{244}\wp_{224}
  -\lambda_2 \wp_{2244}\wp_{24}
  +\tfrac12\lambda_2 \lambda_4 \wp_{2233}\\
  &\quad-\wp_{1223}\wp_{44}\wp_{24}
  -\wp_{13}\wp_{44}\wp_{2244}-2\wp_{13}\wp_{2244}\wp_{24}-
  \wp_{23}\wp_{33}\wp_{2223} \\
&\quad-\tfrac76
  \wp_{33}\wp_{14}\wp_{2233}-\tfrac23 \wp_{33}\wp_{124}\wp_{233}
  -\tfrac43\wp_{33}\wp_{234}\wp_{123}+\tfrac13
  \wp_{34}\wp_{33}\wp_{1233} \\
&\quad+2 \wp_{34}\wp_{233}\wp_{123}+\tfrac13
  \wp_{34}\wp_{2233}\wp_{13} +\lambda_4 \wp_{33}\wp_{1233}+2 \lambda_4
  \wp_{233}\wp_{123} \\
&\quad+\lambda_4 \wp_{2233}\wp_{13} -\lambda_3
  \wp_{233}\wp_{223} -\tfrac34 \lambda_3 \wp_{2223}\wp_{33}
  - \tfrac38 \wp_{1224}\wp_{22}+\tfrac12 \wp_{1222}\wp_{24}\\
  &\quad-\tfrac12 \wp_{1122}\wp_{44}+\tfrac12 \wp_{12}\wp_{2224}
  -\tfrac12 \wp_{11}\wp_{2244}+\tfrac12 \wp_{2222}\wp_{22} \\
  &\quad+\tfrac{1}{12}\lambda_3 \wp_{222333} -\tfrac38
  \wp_{14}\wp_{2222}-\tfrac16 \lambda_4 \wp_{122333} -\tfrac16
  \wp_{34}\wp_{122333}.
\end{aligned}
\end{equation}
\end{cor}
We expect that with better understanding of the PDEs satisfied by the
$\wp_{ij}$, the right hand side of this formula could be simplified
considerably.
\end{remark}

\section{Towards a three-term addition theorem} 

There must be a second main type of addition result satisfied by the
$\sigma$-function of any cyclic trigonal curve:
\begin{remark} 
\label{3term}
The following function associated with {\rm (\ref{cyclic35})},{\rm :}
\begin{equation*}
  \frac{\sigma(u+v+w)\sigma(u+[\zeta]v+[\zeta^2]w)
    \sigma(u+[\zeta^2]v+[\zeta]w)}
  {\sigma(u)^3\sigma(v)^3\sigma(w)^3}
\end{equation*}
is Abelian in $u$, $v$ and $w$, for it is a quotient of third-order
theta-functions in each argument.
\end{remark}
It has triple poles where any of the arguments $u$, $v$ or $w$ $\in$
$\Theta_3$. The left-hand side thus belongs to
$$
   \varGamma(J\times J\times J, \mathcal{O}(3((\Theta^{[2]}\times J
   \times J)\cup(J\times \Theta^{[2]} \times J)\cup (J\times J \times
   \Theta^{[2]}))).
$$ 
It must therefore possess an expansion of the form 
\begin{equation}
\begin{aligned}
  &\frac{\sigma(u+v+w)\sigma(u+[\zeta]v+[\zeta^2]w)
    \sigma(u+[\zeta^2]v+[\zeta]w)}
  {\sigma(u)^3\sigma(v)^3\sigma(w)^3}\\
  &\qquad\qquad= \sum_{i=1}^{81} \sum_{j=1}^{81} \sum_{k=1}^{81}c_{ijk}
  U_i(u)V_j(v)W_k(w).
\end{aligned}\label{3term3}
\end{equation}
where the functions $U_i$, $V_j$ $W_k$ are all basis functions for the 
space  $\varGamma(J, \mathcal{O}(3 \Theta^{[2]}))$.
\begin{remark}
This calculation may be reduced in complexity, by first finding an analogous,
but less symmetric, expansion of the form
\begin{equation}
\begin{aligned}
   &\frac{\sigma(u+v)\sigma(u+[\zeta]v)\sigma(u+[\zeta^2]v)}
{\sigma(u)^3\sigma(v)^3}\\
&\qquad\qquad= \sum_{i=1}^{81} \sum_{j=1}^{81} d_{ij}  U_i(u)V_j(v),
\end{aligned}\label{3term2}
\end{equation}
we note that the Sato weight of the left hand side of (\ref{3term2})
is $-24$, so for each term in the sum, the Sato weights of the
coefficient $d_{ij}$ and the basis functions $U_i(u)$, $V_j(v)$ must
sum to $-24$ also.
\end{remark}
\begin{remark} 
By applying 
\begin{equation}
   \frac13
   \frac{\partial}{\partial u_i}
   \Big(\frac{\partial}{\partial u_j}
  +\frac{\partial}{\partial v_j}
  +\frac{\partial}{\partial w_j}
   \Big)\log
\end{equation}
to a formula like (\ref{3term}), we would obtain algebraic addition
formulae of another type, including a triple-angle formula.  It would
be interesting to compare such formulae with those of (\ref{R7.4}).
\end{remark}
Again, the difficulty of constructing the three-term formula explicitly 
should be considerably reduced when we note that every term in the sum 
must have total Sato weight $-48$, and that both sides are symmetric in $u$,
$v$, and $w$, and invariant under $u \mapsto [\zeta]u$.

Such 3-term addition formulae have been found explicitly for the cases of
the equianharmonic elliptic curve and for the cyclic $(3,4)$ curve; this
can be done wherever a basis of $\varGamma(J, \mathcal{O}(3 \Theta^{[2]}))$
is known. As the dimension of this space is $3^g$, such calculations for
higher genus curves will rapidly become very unwieldy.

\section*{Acknowledgements}
We would like to thank Prof V.Z. Enolski, Dr. S. Matsutani and Prof.
E. Previato for many valuable discussions about this and related work.
Sadie Baldwin was supported by a DTA studentship from EPSRC. Chris
Eilbeck began working on this problem while attending the programme in
Nonlinear Waves at the Mittag-Leffler Institute in Stockholm in 2005,
and he would like to thank Professor H.\ Holden of Trondheim and the
Swedish Academy of Sciences for making this possible.  Some of the
calculations described in this paper were carried out using
Distributed Maple \cite{smb03}, and we are grateful to the author of
this package, Professor Wolfgang Schreiner of RISC-Linz, for help and
advice.


\end{document}